\documentclass[a4paper,12pt]{article}
\usepackage{comment,amsmath,amssymb,amsthm,pstricks,changepage,pifont,graphicx,hyperref,soul,xcolor}
\usepackage[english]{babel}

\setstcolor{red}

\newtheorem{theorem}{Theorem}

\newtheorem{lemma}[theorem]{Lemma}

\newtheorem{remark}[theorem]{Remark}
\newtheorem{definition}[theorem]{Definition}

\newtheorem{example}[theorem]{Example}

\newcommand{\conv}{\text{conv}}
\newcommand{\Tor}{\text{Tor}}
\newcommand{\lw}{\text{lw}}

\begin{document}
\title{Intrinsicness of the Newton polygon for smooth curves on $\mathbb{P}^1\times \mathbb{P}^1$}
\author{Wouter Castryck and Filip Cools}
\date{}

\maketitle

\noindent \begin{abstract} 
\noindent Let $C$ be a smooth projective curve in $\mathbb{P}^1\times \mathbb{P}^1$ of genus $g\neq 4$, and assume that it is birationally equivalent to a curve defined by a Laurent polynomial that is non-degenerate with respect to its Newton polygon $\Delta$. Then we show that the convex hull $\Delta^{(1)}$ of the interior lattice points of $\Delta$ is a standard rectangle, up to a unimodular transformation. Our main auxiliary result, which we believe to be interesting in its own right, is that the first scrollar Betti numbers of $\Delta$-non-degenerate curves are encoded in the combinatorics of $\Delta^{(1)}$, if $\Delta$ satisfies some mild conditions. \\

\end{abstract}

\section{Introduction}

Let $f\in k[x^{\pm 1},y^{\pm 1}]$ be an irreducible Laurent polynomial over an algebraically closed field $k$ of characteristic zero and let $U(f)$ be the curve it defines in the two-dimensional torus $\mathbb{T}^2=(k^*)^2$. The \emph{Newton polygon} $\Delta=\Delta(f)$ of $f$ is the convex hull in $\mathbb{R}^2$ of all the exponent vectors in $\mathbb{Z}^2$ of the monomials that appear in $f$ with a non-zero coefficient. We will always assume that $\Delta$ is two-dimensional. We say that $f$ is \emph{non-degenerate} with respect to its Newton polygon $\Delta$ (or more briefly, $f$ is $\Delta$-non-degenerate) if and only if for each face $\tau\subset\Delta$ (including $\tau=\Delta$) the system 
$$f_\tau=\frac{\partial f_\tau}{\partial x}=\frac{\partial f_\tau}{\partial y}=0$$ does not have any solutions in $\mathbb{T}^2$. Here, $f_\tau$ is obtained from $f$ by only considering the terms that are supported on $\tau$. This condition is generically satisfied. 
Consider the map $$\varphi_\Delta: \mathbb{T}^2 \to \mathbb{P}^{\sharp(\Delta\cap\mathbb{Z}^2)-1} : (x,y) \mapsto (x^i y^j)_{(i,j)\in\Delta\cap\mathbb{Z}^2}.$$ The Zariski closure of its full image $\varphi_\Delta(\mathbb{T}^2)$ is a \emph{toric surface} $\Tor(\Delta)$, while the Zariski closure of $\varphi_\Delta(U(f))$ is a hyperplane section $C_f$ of $\Tor(\Delta)$, which is smooth if $f$ is non-degenerate. We will denote the projective coordinates of $\mathbb{P}^{\sharp(\Delta\cap\mathbb{Z}^2)-1}$ by $X_{i,j}$ where $(i,j)$ runs over $\Delta\cap\mathbb{Z}^2$.

We say that a smooth curve $C$ is $\Delta$-\emph{non-degenerate} if and only if it is birationally equivalent to $U(f)$ for a $\Delta$-non-degenerate Laurent polynomial $f$. Note that if $C$ is moreover projective, then it is isomorphic to $C_f$. If $C$ is $\Delta$-non-degenerate, then a lot of its geometric properties are encoded in the combinatorics of the lattice polygon $\Delta$. For instance, its geometric genus $g(C)$ equals the number of interior lattice points of $\Delta$ \cite{Khov}. Similar interpretations were recently provided for the gonality \cite{CaCo1,Kawaguchi}, the Clifford index and dimension \cite{CaCo1,Kawaguchi}, the scrollar invariants associated to a gonality pencil \cite{CaCo1} and Schreyer's tetragonal invariants \cite{CaCo3}. 

Given this long list, the following question (initiated in \cite{CaCo3}) naturally arises: to what extent can we recover $\Delta$ from the geometry of a $\Delta$-non-degenerate curve? At least, we have to allow two relaxations to this question. First, we can only expect to find back the polygon $\Delta$ up to a \emph{unimodular transformation}, i.e. an affine map of the form 
$$\chi : \mathbb{R}^2\to\mathbb{R}^2 : \begin{pmatrix} x \\ y \end{pmatrix} \mapsto A \begin{pmatrix} x \\ y \end{pmatrix} + B$$
with $A\in \text{GL}_2(\mathbb{Z})$ and $B\in\mathbb{Z}^2$, since these maps correspond to automorphisms of $\mathbb{T}^2$. Secondly, we can (usually) only hope to recover the convex hull of the interior lattice points of $\Delta$, denoted by $\Delta^{(1)}$ (see \cite{CaCo3} for an easy example demonstrating the need for this relaxation). In fact, all the aforementioned combinatorial interpretations are in terms of the combinatorics of $\Delta^{(1)}$ rather than $\Delta$ (e.g. $g(C)=\sharp(\Delta^{(1)}\cap\mathbb{Z}^2)$). 

Given a $\Delta$-non-degenerate curve $C$, we say that the Newton polygon $\Delta$ is \emph{intrinsic} for $C$ if and only if for all $\Delta'$-non-degenerate curves $C'$ that are birationally equivalent to $C$, we have that $\Delta^{(1)}\cong \Delta'^{(1)}$. Hereby, we use $\cong$ to denote the unimodular equivalence relation. 
Before stating some intrinsicness results, we give notations for some special lattice polygons:
\begin{align*}
\square_{a,b} &= \conv\{(0,0),(a,0),(0,b),(a,b)\} \text{ for } a,b\in\mathbb{Z}_{\geq 0},\\
\Sigma &= \conv\{(0,0),(1,0),(0,1)\}, \\
\Upsilon &= \conv\{(-1,-1),(1,0),(0,1)\}.
\end{align*}
The Newton polygon is intrinsic for all rational ($\Delta^{(1)}=\emptyset$), hyperelliptic ($\Delta^{(1)}$ is one-dimensional, and therefore determined by the genus) and trigonal curves of genus at least $5$ ($\Delta^{(1)}$ has lattice width $1$, and is determined by the Maroni invariants). However, there are trigonal curves of genus $4$ for which $\Delta$ is not intrinsic: there exist curves which are non-degenerate with respect to polygons $\Delta$ and $\Delta'$, with $\Delta^{(1)}=\Upsilon$ and $\Delta'^{(1)}=\square_{1,1}$. Intrinsicness of the Newton polygon for tetragonal curves was studied in \cite{CaCo3}: the Newton polygon $\Delta$ is intrinsic if $g(C)\bmod 4\in\{2,3\}$, but it might occasionally be not intrinsic if $g(C)\bmod 4\in\{0,1\}$. From \cite{CaCo1}, it follows that non-degenerate smooth plane curves of degree $d\geq 3$ ($\Delta^{(1)}\cong (d-3)\Sigma$) and curves with Clifford dimension $3$ ($\Delta^{(1)}\cong 2\Upsilon$) have an intrinsic Newton polygon. Moreover, a partial result was given for non-degenerate curves on Hirzebruch surfaces $\mathcal{H}_n$: the value $n$ is intrinsic. 

In this paper, we examine intrinsicness of $\Delta$ for curves on $\mathbb{P}^1\times \mathbb{P}^1$. Namely, we will show that a $\Delta$-non-degenerate curve $C$ of genus $g\neq 4$ can be embedded in $\mathbb{P}^1\times \mathbb{P}^1$ (if and) only if $\Delta^{(1)}=\emptyset$ or $\Delta^{(1)}\cong\square_{a,b}$ for $a,b\in\mathbb{Z}_{\geq 0}$ satisfying $g=(a+1)(b+1)$; see Theorem \ref{P1xP1theorem} in Section \ref{sec_intrinsicness}. In order to prove this result, we give a combinatorial interpretation for the first scrollar Betti numbers of $\Delta$-non-degenerate curves with respect to a gonality pencil, as soon as $\Delta$ satisfies some mild conditions (see Section \ref{sec2}). 

\textbf{Notations.} Let $\mathbb{P}^N$ be a projective space with coordinates $(X_0:\ldots:X_N)$. For each projective variety $V\subset \mathbb{P}^N$, we write $\mathcal{I}(V)\subset k[X_0,\ldots,X_N]$ to indicate the homogeneous ideal of $V$ and $\mathcal{I}_d(V)\subset \mathcal{I}(V)$ to indicate its homogeneous degree $d$ piece. If $J\subset k[X_0,\ldots,X_N]$ is a homogeneous ideal, then $\mathcal{Z}(J)\subset \mathbb{P}^N$ is the zero locus of the polynomials in $J$.

\textbf{Acknowledgements.} We would like to thank Christian Bopp, Marc Coppens and Jeroen Demeyer for some interesting discussions. This research was conducted in the framework of Research Project G093913N of the Research Foundation - Flanders (FWO).

\section{First scrollar Betti numbers} \label{sec2}

\subsection{Definition} \label{subsec2.1}

We start by recalling the definition and some properties of \emph{rational normal scrolls}. 

Let $n \in \mathbb{Z}_{\geq 2}$ and let $\mathcal{E}=\mathcal{O}(e_1)\oplus \cdots \oplus \mathcal{O}(e_n)$ be a locally free sheaf of rank $n$ on $\mathbb{P}^1$. Denote by $\pi:\mathbb{P}(\mathcal{E})\to\mathbb{P}^1$ the corresponding $\mathbb{P}^{n-1}$-bundle. We assume that 
$0\leq e_1\leq e_2\leq \ldots\leq e_n$ and that $e_1 + e_2 + \dots + e_n \geq 2$. Set $N=e_1+e_2+\ldots+e_n+n-1$. Then the image 
$S=S(e_1,\ldots,e_n)$ of the induced morphism
$$\mu:\mathbb{P}(\mathcal{E})\to \mathbb{P}H^0(\mathbb{P}(\mathcal{E}),\mathcal{O}_{\mathbb{P}(\mathcal{E})}(1)),$$
when composed with an isomorphism $\mathbb{P}H^0(\mathbb{P}(\mathcal{E}),\mathcal{O}_{\mathbb{P}(\mathcal{E})}(1)) \rightarrow \mathbb{P}^N$, 
is called a rational normal scroll of type $(e_1, \dots, e_n)$. Up to automorphisms of $\mathbb{P}^N$, rational normal scrolls are fully characterized by their type.

The scroll is smooth if and only if $e_1>0$. In this case, $\mu : \mathbb{P}(\mathcal{E}) \rightarrow S$ is an isomorphism.
If $0=e_1=\ldots=e_\ell<e_{\ell+1}$ with $1 \leq \ell<n$, then the scroll is a cone with an $(\ell-1)$-dimensional vertex.
In this case $\mu : \mathbb{P}(\mathcal{E}) \rightarrow S$ is a resolution of singularities and 
$$\mu_\lambda : \mathbb{P}(\mathcal{E})\cong \mathbb{P}(\mathcal{E}\otimes\mathcal{O}_{\mathbb{P}^1}(\lambda)) \to S'=S(e_1+\lambda,\ldots,e_n+\lambda)$$ is an isomorphism for all integers $\lambda>0$.

The Picard group of $\mathbb{P}(\mathcal{E})$ is freely generated by the class $H$ of a
hyperplane section (more precisely, the class corresponding to $\mu^\ast \mathcal{O}_{\mathbb{P}^N}(1)$) and
the class $R$ of a fiber of $\pi$; i.e. $$\text{Pic}(\mathbb{P}(\mathcal{E}))=\mathbb{Z}H \oplus \mathbb{Z}R.$$
We have the following intersection products: $$H^n=e_1+\ldots+e_n,\ H^{n-1}R=1 \ \text{ and }\ R^2=0$$
(where $R^2 = 0$ means that any appearance of $R^2$ annihilates the product). If we denote the class which corresponds to 
$\mu_\lambda^\ast \mathcal{O}_{\mathbb{P}^{N+n\lambda}}(1)$ by $H'$, we obtain the equality $H'=H+\lambda R$ in $\text{Pic}(\mathbb{P}(\mathcal{E}))$. \\ 

Let $C / k$ be a smooth projective curve of genus $g$ and gonality $\gamma \geq 4$. Assume that $C$ is canonically embedded in $\mathbb{P}^{g-1}$ and fix a gonality pencil $g^1_\gamma$ on $C$. By \cite[Thm. 2]{EiHa}, $$S=\bigcup_{D\in g^1_{\gamma}}\,\langle D\rangle \subset \mathbb{P}^{g-1}$$ is a $(\gamma-1)$-dimensional rational normal scroll containing $C$. If $S$ is of type $(e_1,\ldots,e_{\gamma-1})$, the numbers $e_1,\ldots,e_{\gamma-1}$ are called the \emph{scrollar invariants} of $C$ with respect to $g_{\gamma}^1$. Using the Riemann-Roch theorem, one can see that $e_{\gamma-1}\leq \frac{2g-2}{\gamma}$. Let $\mu : \mathbb{P}(\mathcal{E}) \rightarrow S$ be the corresponding resolution and let $C'$ be the strict transform under $\mu$ of our canonical model. 

The following theorem extends a result from \cite{Schreyer} on tetragonal and pentagonal curves to arbitrary curves.

\begin{theorem} \label{thm_fsBn}
There exist effective divisors $D_1,\ldots,D_{(\gamma^2-3\gamma)/2}$ on $\mathbb{P}(\mathcal{E})$
along with integers $b_1, \dots, b_{(\gamma^2-3\gamma)/2}$, such that
$C'$ is the (scheme theoretic) intersection of the $D_\ell$'s, and such
that for all $\ell$ we have $D_\ell\sim 2H-b_\ell R$.
Moreover, the multiset $\{b_1,\ldots,b_{(\gamma^2-3\gamma)/2}\}$ does
 not depend on the choice of the $D_\ell$'s, and
 \[ \sum_{\ell=1}^{(\gamma^2-3\gamma)/2} b_\ell \, = \, (\gamma-3)g-(\gamma^2-2\gamma-3).\]
\end{theorem}

\begin{proof}
Define $\beta_i=\frac{i(\gamma-2-i)}{\gamma-1}{\gamma \choose i+1}$ and note that $\beta_1 = (\gamma^2-3\gamma)/2$.
The existence follows from \cite[Cor.~4.4]{Schreyer} and its proof, where
the $D_\ell$'s give rise to an exact sequence
of $\mathcal{O}_{\mathbb{P}(\mathcal{E})}$-modules
\begin{multline} \label{schreyer_resolution}
0\to \mathcal{O}_{\mathbb{P}(\mathcal{E})}(-\gamma H+(g-\gamma+1)R)\to \sum_{\ell=1}^{\beta_{\gamma-3}}\mathcal{O}_{\mathbb{P}(\mathcal{E})}(-(\gamma-2)H+b_\ell^{(\gamma-3)}R)\to\cdots \\ \to \sum_{\ell=1}^{\beta_2}\mathcal{O}_{\mathbb{P}(\mathcal{E})}(-3H+b_\ell^{(2)}R)\to\sum_{\ell=1}^{\beta_1}\mathcal{O}_{\mathbb{P}(\mathcal{E})}(-2H+b_\ell R)\to \mathcal{O}_{\mathbb{P}(\mathcal{E})}\to \mathcal{O}_{C'}\to 0.
\end{multline}
Tensoring (\ref{schreyer_resolution}) with $\mathcal{O}_{\mathbb{P}(\mathcal{E})}(2H + bR)$ for a sufficiently large integer $b$ and
computing the Euler characteristics of the terms in the resulting exact sequence, one can show that
\begin{equation} \label{sumbells}
 \sum_\ell b_\ell \, = \, (\gamma-3)g-(\gamma^2-2\gamma-3);
\end{equation}
see \cite[Prop. 2.9]{BoppHoff}.
To conclude the proof, note that the exact sequence
\[ \sum_{\ell = 1}^{\beta_1} \mathcal{O}_{\mathbb{P}(\mathcal{E})}(-D_\ell) \rightarrow \mathcal{O}_{\mathbb{P}(\mathcal{E})} \rightarrow \mathcal{O}_{C'} \rightarrow 0 \]
can be extended to a minimal free resolution of $C'$ on $\mathbb{P}(\mathcal{E})$ of the form (\ref{schreyer_resolution}) by 
minimally resolving the kernel of $\sum_{\ell = 1}^{\beta_1} \mathcal{O}_{\mathbb{P}(\mathcal{E})}(-D_\ell) \rightarrow \mathcal{O}_{\mathbb{P}(\mathcal{E})}$ 
in terms of $\mathcal{O}_{\mathbb{P}(\mathcal{E})}$-modules
(see \cite[proof of Thm.~3.2]{Schreyer}). This resolution is unique up to isomorphism.
\end{proof}

We call the invariants $b_1,\ldots,b_{(\gamma^2-3\gamma)/2}$ the \emph{first scrollar Betti numbers of $C$ with respect to $g^1_\gamma$}. The main goal of this section is to give a combinatorial interpretation for these invariants for non-degenerate curves.

In \cite{CaCo3}, we already treated the case of tetragonal $\Delta$-non-degenerate curves: the first scrollar Betti numbers are given by 
$$\sharp(\partial\Delta^{(1)}\cap\mathbb{Z}^2)-4 \quad \text{and} \quad \sharp(\Delta^{(2)}\cap\mathbb{Z}^2)-1.$$ These numbers are independent from the choice of the gonality pencil. This will no longer be true for non-degenerate curves of higher gonality.

\subsection{Scrollar invariants for non-degenerate curves} \label{subsec2.2}

Let $f$ be a $\Delta$-non-degenerate Laurent polynomial and consider the corresponding smooth curve $C_f\subset \Tor(\Delta)\subset \mathbb{P}^N$ with $N=\sharp(\Delta\cap\mathbb{Z}^2)-1$. Assume that the polygon $\Delta^{(1)}$ is two-dimensional. 

By \cite{Khov}, $C_f$ is a non-rational and non-hyperelliptic curve and there exists a canonical divisor $K_\Delta$ on $C_f$ such that $$H^0(C_f,K_\Delta)=\langle x^iy^j\rangle_{(i,j)\in\Delta^{(1)}\cap\mathbb{Z}^2}$$ (where $x,y$ are functions on $C_f$ through $\varphi_\Delta$). In particular, the curve $C_f$ has genus $g=\sharp(\Delta^{(1)}\cap\mathbb{Z}^2)\geq 3$; see \cite{CaCo1} for more details. Moreover, the Zariski closure $C=C_f^{can}$ of the image of $U(f)$ under 
\begin{equation} \varphi_{\Delta^{(1)}}:\mathbb{T}^2\hookrightarrow \mathbb{P}^{g-1}:(x,y)\mapsto (x^iy^j)_{(i,j)\in\Delta^{(1)}\cap\mathbb{Z}^2} \label{canemb} \end{equation}
is a canonical model for $C_f$. We end up with the inclusions $$C\subset T=\Tor(\Delta^{(1)})=\overline{\varphi_{\Delta^{(1)}}(\mathbb{T}^2)}\subset \mathbb{P}^{g-1},$$
where $T$ is a toric surface since $\Delta^{(1)}$ is two-dimensional.

A \emph{lattice direction} is a primitive integer vector $v=(a,b)\in \mathbb{Z}^2$. The \emph{width} $w(\Delta,v)$ of $\Delta$ with respect to a lattice direction $v$ is the smallest integer $\ell$ such that $\Delta$ is contained in the strip $k\leq aY-bX\leq k+\ell$ of $\mathbb{R}^2$ for some $k\in \mathbb{Z}$. The \emph{lattice width} is defined as $\lw(\Delta)=\min_v w(\Delta,v)$. Lattice directions $v$ that attain the minimum are called \emph{lattice width directions}. 

In \cite{CaCo1}, we gave a combinatorial interpretation for the gonality $\gamma$ of $C=C_f^{can}$ (or $C_f$) in terms of the lattice width of $\Delta$: 
$$\gamma=\begin{cases} 
\lw(\Delta)=\lw(\Delta^{(1)})+2 & \text{if } \Delta\not\cong2\Upsilon \text{ and } \Delta\not\cong d\Sigma  \text{ for all } d\in\mathbb{Z}_{\geq 4}, \\ 
\lw(\Delta)-1=\lw(\Delta^{(1)})+2 & \text{if } \Delta\cong d\Sigma \text{ for some } d\in\mathbb{Z}_{\geq 4}, \\ 
\lw(\Delta)-1=\lw(\Delta^{(1)})+1 & \text{if } \Delta\cong 2\Upsilon,
\end{cases}$$
where we use our assumption that $\Delta^{(1)}$ is two-dimensional. From now on, we make the stronger assumption that $\gamma=\lw(\Delta)\geq 4$, and that $\Delta^{(1)}$ is not equivalent with $(d-3)\Sigma$ or $\Upsilon$, hence $\Delta\not\cong d\Sigma$ and $\Delta\not\cong 2\Upsilon$. Then each lattice width direction $v=(a,b)$ gives rise to a rational map $$C\dashrightarrow \mathbb{P}^1: (x^iy^j)_{(i,j)\in\Delta^{(1)}\cap\mathbb{Z}^2} \mapsto x^ay^b$$ of degree equal to the gonality $\gamma$. We call the corresponding linear pencil $g_\gamma^1$ of $C$ a \emph{combinatorial gonality pencil}. If $\Delta$ is sufficiently big (for a precise statement, see \cite[Corollary 6.3]{CaCo1}), each gonality pencil on $C$ is combinatorial. 

Fix a lattice width direction $v$ of $\Delta$. After applying a suitable unimodular transformation $\chi$, we may assume that $v=(1,0)$ and that $\Delta$ is contained in the horizontal strip $0\leq Y\leq \gamma$ in $\mathbb{R}^2$. So, the gonality map $C\dashrightarrow \mathbb{P}^1$ associated to $v$ is the vertical projection to the $x$-axis. Write $$i^{(-)}(j)=\min\{i\in\mathbb{Z}\,|\,(i,j)\in \Delta^{(1)}\} \text{ and } i^{(+)}(j)=\max\{i\in\mathbb{Z}\,|\,(i,j)\in \Delta^{(1)}\}$$ for all $j\in\{1,\ldots,\gamma-1\}$.
By \cite[Theorem 9.1]{CaCo1}, the scrollar invariants $e_1,\ldots,e_{\gamma-1}$ of $C$ with respect to $g_\gamma^1$ are equal to 
$E_j := i^{(+)}(j)-i^{(-)}(j)$ for $j\in\{1,\dots,\gamma-1\}$
(up to order). In fact, a Zariski dense part of the scroll $S$ is parametrized by 
$$(a_1,\ldots,a_{\gamma-1},x)\in \mathbb{T}^{\gamma}\mapsto (a_j x^{i-i^{(-)}(j)})_{(i,j)\in\Delta^{(1)}\cap\mathbb{Z}^2}=(a_j,\ldots,a_j x^{E_j})_{1\leq j\leq \gamma-1}\in \mathbb{P}^{g-1}.$$ Note that $T=\Tor(\Delta^{(1)})\subset S$ since the map $\varphi_{\Delta^{(1)}}$ can be obtained from the above parametrization by restricting to $a_j=x^{i^{(-)}(j)}y^j$, so we get the inclusions 
\begin{equation} \label{incl1} C\subset T\subset S\subset \mathbb{P}^{g-1}. \end{equation}

If $S$ is singular, then $\mu: S'=S(e_1+\lambda,\ldots,e_{\gamma-1}+\lambda)\cong \mathbb{P}(\mathcal{E}) \to S$ is a resolution of singularities for each integer $\lambda>0$ (hereby, we slightly abuse notation: $\mu$ is the map $\mu\circ\mu_\lambda^{-1}$ using the notations in Section \ref{subsec2.1}). Let $C'$ and $T'$ be the strict transforms of respectively $C$ and $T$ under $\mu$. 
For each lattice polygon $\Gamma\subset\mathbb{R}^2$, write $\Gamma[\lambda]$ to denote the Minkowski sum of $\Gamma$ and $[(0,0),(\lambda,0)]\subset \mathbb{R}^2$. In other words, $\Gamma[\lambda]$ is obtained from $\Gamma$ by stretching it out in the horizontal direction over a distance $\lambda$. Using this notation, one can see that $T'=\Tor(\Delta^{(1)}[\lambda])=\Tor(\Delta[\lambda]^{(1)})$. We end up with the inclusions \begin{equation} \label{incl2} C'\subset T'\subset S'\subset \mathbb{P}^{g-1+\lambda(\gamma-1)}. \end{equation}

\subsection{First scrollar Betti numbers of toric surfaces} \label{subsec2.3}

Let $C$ be a $\Delta$-non-degenerate curve
and fix a combinatorial gonality pencil $g^1_\gamma$ on $C$, corresponding to a lattice direction $v$. 
We work under the following assumptions:
\begin{enumerate}
\item[(i)] $\Delta^{(1)}$ is not equivalent with $(d-3)\Sigma$ for any $d\geq 3$ or $\Upsilon$, and $\gamma=\lw(\Delta)\geq 4$, 
\item[(ii)] $v=(1,0)$ and $\Delta$ is contained in the horizontal strip $0\leq Y\leq \gamma$,
so that $g^1_\gamma$ corresponds to vertical projection,
\item[(iii)] the curve $C$ is canonically embedded, so that we obtain the sequence of inclusions $C \subset T \subset S \subset \mathbb{P}^{g-1}$ from \eqref{incl1}.
\end{enumerate}
Recall that the scrollar invariants $e_1, \dots, e_{\gamma - 1}$ of $C$ with respect to $g^1_\gamma$ match with $E_1,\ldots,E_{\gamma-1}$ (up to order). 
Consider $\mu : \mathbb{P}(\mathcal{E}) \rightarrow S$ and let $T'\subset \mathbb{P}(\mathcal{E})$ be the strict transform of 
$T=\text{Tor}(\Delta^{(1)})$ under $\mu$, as in Section \ref{subsec2.2}.
If $\Delta$ satisfies the condition $\mathcal{P}_1(v)$ defined below (see Definition \ref{def_P1}), we will provide effective divisors $D_1,\ldots,D_{{\gamma-2 \choose 2}}$ on $\mathbb{P}(\mathcal{E})$
along with integers $b_1, \dots, b_{{\gamma-2 \choose 2}}$, such that the following three conditions are satisfied:
\begin{itemize}
\item $T'$ is the (scheme theoretic) intersection of the $D_\ell$'s, 
\item $D_\ell\sim 2H-b_\ell R$ for all $\ell$, and,
\item $\sum_{\ell=1}^{{\gamma-2 \choose 2}} b_\ell \, = (\gamma-4)g-(\gamma^2-3\gamma)+\sharp(\partial\Delta^{(1)}\cap\mathbb{Z}^2)$.
\end{itemize}
In what follows, we will also assume that $e_1 > 0$, so that $\mathbb{P}(\mathcal{E})\cong S$. This condition
is not essential (see Remark \ref{rmk_ezero}), but it allows us to work with the inclusion $T\subset S$ rather than $T'\subset \mathbb{P}(\mathcal{E})$. 
For convenience, we will denote the invariants by $B_{j_1,j_2}$ and the divisors by $D_{j_1,j_2}$, where $j_1,j_2\in\{1,\ldots,\gamma-1\}$ such that $j_2-j_1\geq 2$. \\

For each $j_1,j_2\in\{1,\ldots,\gamma-1\}$ such that $j_2-j_1\geq 2$ and $1\leq r\leq \frac{j_2-j_1}{2}$, let $Y_{j_1,j_2,r}\subset S$ be the subvariety defined by the binomials of
$\mathcal{I}_2(\text{Tor}(\Delta^{(1)}))$ having the form $$X_{i_1,j_1}X_{i_2,j_2} - X_{i'_1,j_1+r}X_{i'_2,j_2-r}.$$ 
One can see that $Y_{j_1,j_2,r}$ is a $(\gamma-2)$-dimensional toric variety $\text{Tor}(\Omega_{j_1,j_2,r})$, where $\Omega_{j_1,j_2,r}\subset \mathbb{R}^{\gamma-2}$ is a full-dimensional lattice polytope (see Example \ref{exOmega} for a tangible instance). The (Euclidean) volume of this polytope equals
$$\frac{1}{(\gamma-2)!}(2(E_1+\ldots+E_{\gamma-1})-(E_{j_1}+E_{j_2}-\epsilon_{j_1,j_2,r})),$$
where $\epsilon_{j_1,j_2,r}$ is defined as $\epsilon^{(-)}_{j_1,j_2,r}+\epsilon^{(+)}_{j_1,j_2,r}$, with 
\begin{eqnarray*} \label{defepsilons-}
\epsilon^{(-)}_{j_1,j_2,r} &=& 
\begin{cases}
0 \quad \text{ if } & i^{(-)}(j_1+r)+i^{(-)}(j_2-r)\leq i^{(-)}(j_1)+i^{(-)}(j_2) \\
1 \quad \text{ if } & i^{(-)}(j_1+r)+i^{(-)}(j_2-r) > i^{(-)}(j_1)+i^{(-)}(j_2) \\
\end{cases} \\ 
&=& \max\{ 0, (i^{(-)}(j_1+r)+i^{(-)}(j_2-r)) - (i^{(-)}(j_1)+i^{(-)}(j_2)) \}, 
\end{eqnarray*}
and 
\begin{eqnarray*} \label{defepsilons+}
\epsilon^{(+)}_{j_1,j_2,r} &=& 
\begin{cases}
0 \quad \text{ if } & i^{(+)}(j_1+r)+i^{(+)}(j_2-r)\geq i^{(+)}(j_1)+i^{(+)}(j_2)\\
1 \quad \text{ if } & i^{(+)}(j_1+r)+i^{(+)}(j_2-r)< i^{(+)}(j_1)+i^{(+)}(j_2)\\
\end{cases} \\
&=& \max\{ 0, (i^{(+)}(j_1)+i^{(+)}(j_2)) - (i^{(+)}(j_1+r)+i^{(+)}(j_2-r)) \}.
\end{eqnarray*}
In the above equalities for $\epsilon^{(-)}_{j_1,j_2,r}$ and $\epsilon^{(+)}_{j_1,j_2,r}$, we use the following result.

\begin{lemma} \label{lem_interior} The inequalities 
$$i^{(-)}(j_1+r)+i^{(-)}(j_2-r)\leq i^{(-)}(j_1)+i^{(-)}(j_2)+1$$ and
$$i^{(+)}(j_1+r)+i^{(+)}(j_2-r)\geq i^{(+)}(j_1)+i^{(+)}(j_2)-1$$
hold for all $j_1,j_2\in\{1,\ldots,\gamma-1\}$ such that $j_2-j_1\geq 2$ and $1\leq r\leq \frac{j_2-j_1}{2}$.
\end{lemma}

\begin{proof} 
We only show the first inequality; the second one follows by symmetry. Consider the line segment $L=[(i^{(-)}(j_1),j_1),(i^{(-)}(j_2),j_2)]$, and let 
$(i',j_1+r)$ and $(i'',j_2-r)$ be the intersection points of $L$ with the horizontal lines at heights $j_1+r$ and $j_2-r$. 
Note that $L$ is contained in the interior of $\Delta$ and that $i'+i''=i^{(-)}(j_1)+i^{(-)}(j_2)$. If $i^{(-)}(j_1+r)+i^{(-)}(j_2-r)\geq i^{(-)}(j_1)+i^{(-)}(j_2)+2=i'+i''+2$, then $i'\leq i^{(-)}(j_1+r)-1$ or $i''\leq i^{(-)}(j_2-r)-1$, so $(i^{(-)}(j_1+r)-1,j_1+r)$ or $(i^{(-)}(j_2-r)-1,j_2-r)$ is a lattice point lying in the interior of $\Delta$. This is in contradiction with the definition of $i^{(-)}(\cdot)$. 
\end{proof}

\begin{example} \label{exOmega}
Assume that $\Delta = \Delta(f)$ is as in Figure \ref{FigGon5a} (here $\gamma = 5$).
\begin{figure}[h!]
  \centering
    \includegraphics[height=2.5cm]{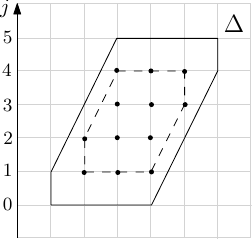}
		\caption{picture of $\Delta$}
		\label{FigGon5a}
\end{figure}

\noindent Appropriate instances of $\Omega_{j_1,j_2,r}$ can be realized as in Figure \ref{FigOmegas}.
\begin{figure}[h!]
  \centering
    \includegraphics[height=2cm]{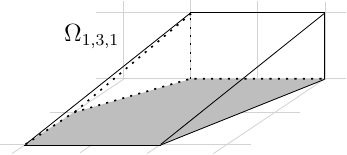} 
		\includegraphics[height=2cm]{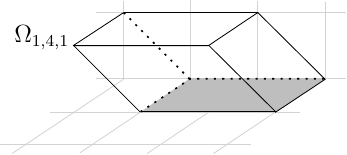} 
		\includegraphics[height=2cm]{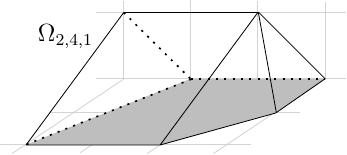} 
		\caption{picture of the $\Omega_{j_1,j_2,r}$'s}
		\label{FigOmegas}
\end{figure}

\noindent Here, $\epsilon_{1,3,1}=1$ (since $\epsilon^{(+)}_{1,3,1}=1$), $\epsilon_{1,4,1}=0$ and $\epsilon_{2,4,1}=1$ (since $\epsilon^{(-)}_{2,4,1}=1$). 
\end{example}

\noindent Since the intersection of $Y_{j_1,j_2,r}$ with a typical fiber of $S\to\mathbb{P}^1$ is a quadratic hypersurface,
we have $Y_{j_1,j_2,r}\sim 2H-B_{j_1,j_2,r}R$ for some $B_{j_1,j_2,r}\in\mathbb{Z}$. Taking the intersection product with $H^{\gamma-2}$, we get
$$ \begin{array}{lll}
\deg Y_{j_1,j_2,r}=Y_{j_1,j_2,r}\cdot H^{\gamma-2}&=&2H^{\gamma-1}-B_{j_1,j_2,r}H^{\gamma-2}R \\ &=&2(e_1+\ldots+e_{\gamma-1})-B_{j_1,j_2,r}
\\ &=& 2(E_1+\ldots+E_{\gamma-1})-B_{j_1,j_2,r},
\end{array}$$
but $\deg Y_{j_1,j_2,r}=(\gamma-2)!\cdot \text{Vol}(\Omega_{j_1,j_2,r})$, so $B_{j_1,j_2,r}$ equals $E_{j_1}+E_{j_2}-\epsilon_{j_1,j_2,r}$. 

Write 
$$\mathcal{S}^{(-)}_{j_1,j_2}=\left\{r\in \left\{1,\ldots,\left\lfloor \frac{j_2-j_1}{2}\right\rfloor\right\}\,|\,\epsilon^{(-)}_{j_1,j_2,r}=0\right\}$$ and 
$$\mathcal{S}^{(+)}_{j_1,j_2}=\left\{r\in \left\{1,\ldots,\left\lfloor \frac{j_2-j_1}{2}\right\rfloor\right\}\,|\,\epsilon^{(+)}_{j_1,j_2,r}=0\right\}.$$

\begin{definition} \label{def_P1}
We say that $\Delta$ satisfies condition $\mathcal{P}_1(v)$ if and only if there are no integers $j_1,j_2\in\{1,\ldots,\gamma-1\}$ with $j_2-j_1\geq 2$ such that $\mathcal{S}^{(-)}_{j_1,j_2}$ and $\mathcal{S}^{(+)}_{j_1,j_2}$ are non-empty and disjoint. 
\end{definition} 

In other words, the condition $\mathcal{P}_1(v)$ means that for each pair of integers $j_1,j_2\in\{1,\ldots,\gamma-1\}$ with $j_2-j_1\geq 2$ either at least one of the sets $\mathcal{S}^{(-)}_{j_1,j_2}$, $\mathcal{S}^{(+)}_{j_1,j_2}$ is empty, or there is a common $r\in \left\{1,\ldots,\left\lfloor \frac{j_2-j_1}{2}\right\rfloor\right\}$ for which $\epsilon^{(-)}_{j_1,j_2,r}=\epsilon^{(+)}_{j_1,j_2,r}=0$.
There is a useful criterion to check whether $\mathcal{S}^{(-)}_{j_1,j_2}$ is empty or not (and analogous for $\mathcal{S}^{(+)}_{j_1,j_2}$): $\mathcal{S}^{(-)}_{j_1,j_2}=\emptyset$ if and only if all the lattice points $(i^{(-)}(j),j)$ with $j_1<j<j_2$ lie strictly right from the line segment $L=\left[(i^{(-)}(j_1),j_1),(i^{(-)}(j_2),j_2)\right]$. 

In the above definition, we also allow the lattice direction $v$ to be different from $(1,0)$: in that case, first take a unimodular transformation $\chi$ such that $\chi(v)=(1,0)$ and that $\chi(\Delta)$ is contained in the horizontal strip $0\leq Y\leq \gamma$, and replace $\Delta$ by $\chi(\Delta)$ while checking the condition. The definition is independent of the particular choice of the unimodular transformation $\chi$.

In fact, in some of the examples below, the lattice direction $v$ is $(1,0)$, but $\Delta$ is contained in a horizontal strip of the form $k\leq Y\leq k+\gamma$ with $k\neq 0$. In that case, we do not really need to apply a unimodular transformation $\chi$ first: we can define the sets $\mathcal{S}^{(-)}_{j_1,j_2}$ and $\mathcal{S}^{(+)}_{j_1,j_2}$ for $j_1,j_2\in\{k+1,\ldots,k+\gamma-1\}$.

\begin{example}
Assume that a part of $\Delta^{(1)}$ looks as in Figure \ref{FigDelta} (for some large enough $n$).
\begin{figure}[h]
  \centering
    \includegraphics[width=0.5\textwidth]{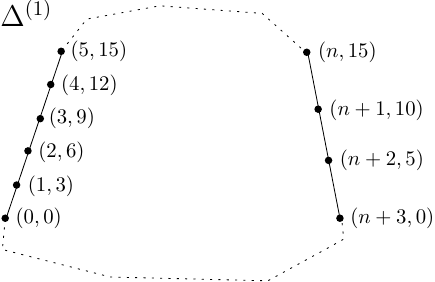}
		\caption{part of $\Delta^{(1)}$}
		\label{FigDelta}
\end{figure}

\noindent In Table \ref{tableS}, the sets $\mathcal{S}^{(-)}_{j_1,j_2}$ and $\mathcal{S}^{(+)}_{j_1,j_2}$ are given for all couples $(j_1,j_2)$ with $j_1+j_2=15$ in this part of the polytope $\Delta^{(1)}$.
\begin{table}[h]
	\centering
  \begin{tabular}{l|c|c}
	$(j_1,j_2)$ & ${\mathcal{S}^{(-)}_{j_1,j_2}}_{}$ & ${\mathcal{S}^{(+)}_{j_1,j_2}}_{}$  \\ \hline \hline 
	$(0,15)$ & $\{3,6\}$ & $\{5\}$ \\ \hline
	$(1,14)$ & $\{1,2,3,4,5,6\}$ & $\{1,2,3,4,5,6\}$ \\ \hline
	$(2,13)$ & $\{1,2,3,4,5\}$ & $\{1,2,3,4,5\}$ \\ \hline
	$(3,12)$ & $\{3\}$ & $\{1,2,3,4\}$ \\ \hline
	$(4,11)$ & $\{1,2,3\}$ & $\{1,2,3\}$ \\ \hline
	$(5,10)$ & $\{1,2\}$ & $\emptyset$ \\ \hline
	$(6,9)$ & $\emptyset$ & $\{1\}$ 
		\end{tabular}
	\caption{table of subsets $\mathcal{S}^{(.)}_{j_1,j_2}$}
	\label{tableS}
\end{table}

\noindent We conclude that $\Delta$ does not satisfy condition $\mathcal{P}_1(v)$ (consider $j_1=0$ and $j_2=15$). 
\end{example}

For all polygons $\Delta$ that satisfy condition $\mathcal{P}_1(v)$, we give a recipe to construct the divisors $D_{j_1,j_2}$ in terms of the subvarieties $Y_{j_1,j_2,r}$.

\begin{definition} \label{defDij}
Assume that the lattice polygon $\Delta$ satisfies condition $\mathcal{P}_1(v)$.
\begin{itemize}
\item If $\mathcal{S}^{(-)}_{j_1,j_2}\cap \mathcal{S}^{(+)}_{j_1,j_2}\neq \emptyset$, we define $D_{j_1,j_2}$ as $Y_{j_1,j_2,r}$ with 
$r\in \mathcal{S}^{(-)}_{j_1,j_2}\cap \mathcal{S}^{(+)}_{j_1,j_2}$ minimal. For instance, if $\epsilon_{j_1,j_2,1}=0$, we have 
$D_{j_1,j_2}=Y_{j_1,j_2,1}$. Set $\epsilon_{j_1,j_2}=\epsilon_{j_1,j_2,r}=0$. 
\item If $\mathcal{S}^{(-)}_{j_1,j_2}=\emptyset$ and $\mathcal{S}^{(+)}_{j_1,j_2}\neq \emptyset$ or vice versa, take 
$r\in \mathcal{S}^{(-)}_{j_1,j_2}\cup \mathcal{S}^{(+)}_{j_1,j_2}$ minimal, define $D_{j_1,j_2}=Y_{j_1,j_2,r}$ and set 
$\epsilon_{j_1,j_2}=\epsilon_{j_1,j_2,r}=1$.
\item If $\mathcal{S}^{(-)}_{j_1,j_2}=\mathcal{S}^{(+)}_{j_1,j_2}=\emptyset$, define $D_{j_1,j_2}=Y_{j_1,j_2,1}$ and set $\epsilon_{j_1,j_2}=\epsilon_{j_1,j_2,1}=2$. 
\end{itemize}
\end{definition}

\begin{remark} \label{rmk_defD}
In Definition \ref{defDij}, the divisor $D_{j_1,j_2}$ is always of the form $Y_{j_1,j_2,r}$ and $r$ is chosen such that 
$\epsilon_{j_1,j_2,r}$ is minimal, or equivalently, $B_{j_1,j_2,r}$ is maximal. Moreover, if $D_{j_1,j_2}=Y_{j_1,j_2,r}$ and if we define $\epsilon_{j_1,j_2}^{(-)}=\epsilon_{j_1,j_2,r}^{(-)}$ and $\epsilon_{j_1,j_2}^{(+)}=\epsilon_{j_1,j_2,r}^{(+)}$, then 
\begin{equation} \label{min-formula}
\epsilon_{j_1,j_2}^{(-)} = \min_s \epsilon_{j_1,j_2,s}^{(-)}, \quad \epsilon_{j_1,j_2}^{(+)}=\min_t \epsilon_{j_1,j_2,t}^{(+)} \quad \text{ and } \quad 
\epsilon_{j_1,j_2} = \epsilon_{j_1,j_2}^{(-)} + \epsilon_{j_1,j_2}^{(+)}.
\end{equation} 
Here, it is crucial that $\Delta$ satisfies condition $\mathcal{P}_1(v)$: if $\mathcal{S}^{(-)}_{j_1,j_2}$ and $\mathcal{S}^{(+)}_{j_1,j_2}$ would be non-empty and disjoint, then $\min_r \epsilon_{j_1,j_2,r} = 1$ (take $r\in \mathcal{S}^{(-)}_{j_1,j_2}\cup \mathcal{S}^{(+)}_{j_1,j_2}$), but $\min_s \epsilon_{j_1,j_2,s}^{(-)}=\min_t \epsilon_{j_1,j_2,t}^{(+)}=0$.

If we set $B_{j_1,j_2}=E_{j_1}+E_{j_2}-\epsilon_{j_1,j_2}$, we have that $D_{j_1,j_2}\sim 2H-B_{j_1,j_2}R$ and  
$$\begin{array}{lll} 
\sum_{j_2-j_1\geq 2}\,B_{j_1,j_2} 
&=&(\gamma-4)(E_1+\ldots+E_{\gamma-1})+E_1+E_{\gamma-1}-\sum_{j_2-j_1\geq 2}\,\epsilon_{j_1,j_2}\\
&=&(\gamma-4)(g-\gamma+1)+E_1+E_{\gamma-1}-\sum_{j_2-j_1\geq 2}\,\epsilon_{j_1,j_2}. 
\end{array}$$
\end{remark}

\begin{example} \label{easyexample}
If $\partial\Delta^{(1)}$ meets each horizontal line on height $j\in\{2,\ldots,\gamma-2\}$ in two lattice points, we have $\epsilon_{j_1,j_2,r}=0$ and $\mathcal{S}^{(-)}_{j_1,j_2}=\mathcal{S}^{(+)}_{j_1,j_2}=\left\{1,\ldots,\left\lfloor \frac{j_2-j_1}{2}\right\rfloor\right\}$ for all $j_1,j_2,r$. Hence, $\Delta$ satisfies condition $\mathcal{P}_1(v)$. Moreover, $\epsilon_{j_1,j_2}=0$ and $D_{j_1,j_2}=Y_{j_1,j_2,1}$ for all $j_1,j_2$.  
In this case, $$\sharp(\partial\Delta^{(1)}\cap\mathbb{Z}^2)=(E_1+1)+(E_{\gamma-1}+1)+2(\gamma-3)$$ and $\sum \epsilon_{j_1,j_2}=0$, so $$\sum B_{j_1,j_2}=(\gamma-4)g-(\gamma^2-3\gamma)+\sharp(\partial\Delta^{(1)}\cap\mathbb{Z}^2).$$ 
\end{example}

\begin{remark} \label{rmk_ezero}
If $S$ is singular, let $\lambda>0$ be an integer and consider the inclusions from \eqref{incl2}. Note that $\Delta[\lambda]$ satisfies condition $\mathcal{P}_1(v)$ if and only if $\Delta$ satisfies condition $\mathcal{P}_1(v)$.
We can define the subvarieties $Y_{j_1,j_2,r}$ and $D_{j_1,j_2}$ of $S'$ in the same way as we did before (using $\Delta[\lambda]$ instead of $\Delta$). Since $H'=H+\lambda R$, we get that $$Y_{j_1,j_2,r}\sim 2H'-((E_{j_1}+\lambda)+(E_{j_2}+\lambda)-\epsilon_{j_1,j_2,r})R=2H-B_{j_1,j_2,r}R$$ and $D_{j_1,j_2}\sim 2H-B_{j_1,j_2}R$. 
\end{remark}

We are now able to state and prove the main result of this subsection.  
 
\begin{theorem} \label{thm_sfBn_toric}
If $\Delta$ satisfies condition $\mathcal{P}_1(v)$, there exist ${\gamma-2 \choose 2}$ effective divisors $D_{j_1,j_2}$ on $\mathbb{P}(\mathcal{E})$ (with $j_1,j_2\in\{1,\ldots,\gamma-1\}$ and $j_2-j_1\geq 2$) such that
\begin{itemize}                                                                               
\item $T'$ is the (scheme theoretic) intersection of the divisors $D_{j_1,j_2}$,
\item $D_{j_1,j_2}\sim 2H-B_{j_1,j_2}R$ for all $j_1,j_2$, where $B_{j_1,j_2} = E_{j_1} + E_{j_2} - \epsilon_{j_1,j_2}$, and, 
\item $\sum_{j_2-j_1\geq 2}\,B_{j_1,j_2}=(\gamma-4)g-(\gamma^2-3\gamma)+\sharp(\partial\Delta^{(1)}\cap\mathbb{Z}^2)$.
\end{itemize}
\end{theorem}

\begin{proof} 
By Remark \ref{rmk_ezero}, we may assume that $S$ is smooth, hence $\mathbb{P}(\mathcal{E})\cong S$. 
We need to prove that $\mathcal{I}(\text{Tor}(\Delta^{(1)}))=\mathcal{I}(\bigcap D_{j_1,j_2})$, where the inclusion $\mathcal{I}(\bigcap D_{j_1,j_2})\subset \mathcal{I}(\text{Tor}(\Delta^{(1)}))$ is trivial. Pick an arbitrary quadratic binomial 
$$f=X_{i_1,j_1}X_{i_2,j_2}-X_{i_3,j_3}X_{i_4,j_4}\in \mathcal{I}(\text{Tor}(\Delta^{(1)})).$$ These binomials generate the ideal, so we only need to show that $f\in \mathcal{I}(\bigcap D_{j_1,j_2})$. Note that $j_1+j_2=j_3+j_4$, so we may assume that $j_1\leq j_3\leq j_4\leq j_2$.
Moreover, if $j_1=j_3$ and $j_4=j_2$, we get that $f\in \mathcal{I}(S)\subset \mathcal{I}(\bigcap D_{j_1,j_2})$. So we may even assume that $j_1<j_3$. 

Take $r$ such that $D_{j_1,j_2}=Y_{j_1,j_2,r}$. We claim that $$I:=i_1+i_2=i_3+i_4\geq i^{(-)}(j_1+r)+i^{(-)}(j_2-r).$$ If $\epsilon_{j_1,j_2}^{(-)}=\epsilon_{j_1,j_2,r}^{(-)}=0$, this follows from $$I\geq i^{(-)}(j_1)+i^{(-)}(j_2)\geq i^{(-)}(j_1+r)+i^{(-)}(j_2-r).$$ If $\epsilon_{j_1,j_2}^{(-)}=\epsilon_{j_1,j_2,r}^{(-)}=1$, we have that $\epsilon_{j_1,j_2,j_3-j_1}^{(-)}=1$ by \eqref{min-formula} (since $\Delta$ satisfies condition $\mathcal{P}_1(v)$), hence $$I\geq i^{(-)}(j_3)+i^{(-)}(j_4)=i^{(-)}(j_1+r)+i^{(-)}(j_2-r),$$ where we use Lemma \ref{lem_interior}. Analogously, we can show that $I\leq i^{(+)}(j_1+r)+i^{(+)}(j_2-r)$.

The above claim implies that we can find integers $i'_1,i'_2$ such that $i'_1+i'_2=I$, $i^{(-)}(j_1+r)\leq i'_1 \leq i^{(+)}(j_1+r)$, $i^{(-)}(j_2-r)\leq i'_2 \leq i^{(+)}(j_2-r)$, hence 
$$X_{i_1,j_1}X_{i_2,j_2}-X_{i'_1,j_1+r}X_{i'_2,j_2-r}\in \mathcal{I}(D_{j_1,j_2})=\mathcal{I}(Y_{j_1,j_2,r}).$$ 
So we may replace the term $X_{i_1,j_1}X_{i_2,j_2}$ in $f$ by $X_{i'_1,j_1+r}X_{i'_2,j_2-r}$ (and in particular, $j_1$ by $j_1+r$ and $j_2$ by $j_2-r$). Continuing in this way, we will eventually get that $j_1=j_3$ and $j_4=j_2$, hence $f\in \mathcal{I}(S)$. This will happen after a finite number of steps since the maximum of $j_2-j_1$ and $j_4-j_3$ decreases after each step. 

We are left with proving the formula for the sum of the $B_{j_1,j_2}$'s. By Remark \ref{rmk_defD} and the elaboration of Example \ref{easyexample}, it suffices to show that the sum of the $\epsilon_{j_1,j_2}$ counts the number of times that 
$\partial\Delta^{(1)}$ intersects the horizontal lines on height $2,\ldots,\gamma-2$ in a non-lattice point. Let $A^{(-)}$ be the set of couples $(j_1,j_2)$ such that $j_1,j_2\in\{1,\ldots,\gamma-1\}$, $j_2-j_1\geq 2$ and $\mathcal{S}_{j_1,j_2}^{(-)}=\emptyset$ (or equivalently, the line segment $L=\left[(i^{(-)}(j_1),j_1),(i^{(-)}(j_2),j_2)\right]$ passes left from all the lattice points $(i^{(-)}(j'),j')$ with $j_1<j'<j_2$). Let $B^{(-)}$ be the set of integers $j\in\{1,\ldots,\gamma-1\}$ such that $(i^{(-)}(j),j)\not \in \partial\Delta^{(1)}$. We claim that the sets $A^{(-)}$ and $B^{(-)}$ have the same cardinality. We will do this by giving a concrete bijection between these sets. Analogously, we can define the sets $A^{(+)}$ and $B^{(+)}$, and prove that they have the same number of elements. The theorem follows directly, since $\sharp(A^{(-)}\cup A^{(+)})=\sum_{j_1,j_2}\,\epsilon_{j_1,j_2}$ by \eqref{min-formula} and $\sharp(B^{(-)}\cup B^{(+)})$ is the number of non-lattice point intersections. 

If $(j_1,j_2)\in A^{(-)}$, then the line segment $L=[(i^{(-)}(j_1),j_1),(i^{(-)}(j_2),j_2)]$ will pass at the left hand side of the lattice points $(i^{(-)}(j),j)$ with $j_1<j<j_2$. For precisely one of these lattice points, the horizontal distance to $L$ will be equal to the minimal value $\frac{1}{j_2-j_1}$. Consider the map 
$$\alpha^{(-)}:A^{(-)}\to B^{(-)}$$ sending the couple $(j_1,j_2)$ to the value of $j$ of that lattice point; see below for an example. On the other hand, if $j\in B^{(-)}$, thus $(i^{(-)}(j),j)\not \in \partial\Delta^{(1)}$, then there should be lattice points $(i^{(-)}(j_1),j_1)$ and $(i^{(-)}(j_2),j_2)$ with $j_1<j<j_2$ such that $L=[(i^{(-)}(j_1),j_1),(i^{(-)}(j_2),j_2)]$ passes left from 
$(i^{(-)}(j),j)$. If we take a couple $(j_1,j_2)$ that satisfies this property and has a minimal value for $j_2-j_1$, then $(j_1,j_2)\in A^{(-)}$. Indeed, if $(i^{(-)}(j'),j')$ with $j_1<j'<j_2$ would lie on or left from $L$, then either $(j_1,j')$ or $(j',j_2)$ would also satisfy the condition and would have a smaller value for the difference of the heights. Now let's show that the couple $(j_1,j_2)$ is unique. If not, there exists another couple $(j_1',j_2')\in A^{(-)}$ with $j_1'<j<j_2'$ such that $L'=[(i^{(-)}(j_1'),j_1'),(i^{(-)}(j_2'),j_2')]$ passes left from $(i^{(-)}(j),j)$ with $j_2'-j_1'=j_2-j_1$. We may assume that $j_1'<j_1<j<j_2'<j_2$. Then $L$ passes left from $(i^{(-)}(j_2'),j_2')$ and $L'$ passes left from $(i^{(-)}(j_1),j_1)$, so $L''=[(i^{(-)}(j_1'),j_1'),(i^{(-)}(j_2),j_2)]$ passes left from all the lattice points $(i^{(-)}(j'),j')$ with $j_1'<j'<j_2$ (see Figure \ref{uniqueness}). 
\begin{figure}[h!]
  \centering
    \includegraphics[height=4cm]{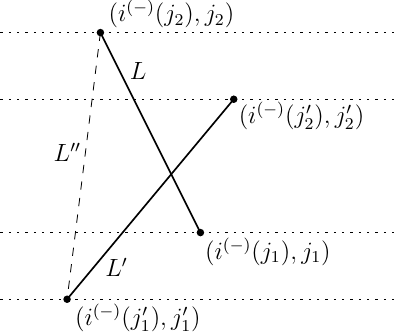}
		\caption{the line segments $L$, $L'$ and $L''$}
		\label{uniqueness}
\end{figure}
Let's denote the horizontal distance from the line segment $L''$ to the lattice point $(i^{(-)}(j'),j')$ by $d(j')$. Using Lemma \ref{lem_interior}, we obtain that $d(j_1'+r)+d(j_2-r)=1$ for all $1\leq r\leq \frac{j_2-j_1'}{2}$. For precisely one integer $j_1'<j'<j_2$, the distance $d(j')$ is equal to the minimal value $\frac{1}{j_2-j_1'}$. On the other hand, except from $j'=j_1$ or $j'=j_2'$, the distance $d(j')$ has be at least $\frac{1}{j_2-j_1}=\frac{1}{j_2'-j_1'}$, since $(i^{(-)}(j'),j')$ lies 
strictly right from $L$ or $L'$. So we may assume that $j'=j_1$ (the case $j'=j_2'$ is analogous), hence 
$d(j_1)=\frac{1}{j_2-j_1'}$ and $d(j_2')=1-\frac{1}{j_2-j_1'}$ (using $r=j_1-j_1'$ above). It follows that the horizontal distance to $L''$ from the point on $L'$ on height $j_1$ is equal to $$\frac{j_1-j_1'}{j_2'-j_1'}\cdot d(j_2')=\frac{j_1-j_1'}{j_2'-j_1'}\cdot \frac{j_2-j_1'-1}{j_2-j_1'}\geq \frac{j_1-j_1'}{j_2-j_1'} \geq \frac{1}{j_2-j_1'}=d(j_1),$$ so $L'$ does not pass left from $(i^{(-)}(j_1),j_1)$, a contradiction. In conclusion we can consider the map $$\beta^{(-)}:B^{(-)}\to A^{(-)}$$ sending $j$ to the unique such couple $(j_1,j_2)$.   

The maps $\alpha^{(-)}$ and $\beta^{(-)}$ are each others inverse. For instance, to prove that the map $\alpha^{(-)}\circ\beta^{(-)}$ is the identity map on $B^{(-)}$, consider $j\in B^{(-)}$ and write $\beta^{(-)}(j)=(j_1,j_2)$. If $\alpha^{(-)}(j_1,j_2)=j'\neq j$, then the horizontal distance from $(i^{(-)}(j),j)$ to $L=[(i^{(-)}(j_1),j_1),(i^{(-)}(j_2),j_2)]$ is of the form $\frac{d}{j_2-j_1}$ with $1<d<j_2-j_1$. But then either $L'=[(i^{(-)}(j'),j'),(i^{(-)}(j_2),j_2)]$ (the case $j'<j$) or $L''=[(i^{(-)}(j_1),j_1),(i^{(-)}(j'),j')]$ (the case $j'>j$) passes left from $(i^{(-)}(j),j)$. This is in contradiction with $\beta^{(-)}(j)=(j_1,j_2)$, since $j_2-j'$ and $j'-j_1$ are both strictly smaller than $j_2-j_1$. We leave the proof of the equality $\beta^{(-)}\circ\alpha^{(-)}=\text{Id}_{A^{(-)}}$ as an exercise. 
\end{proof}

\begin{example}
Consider a polygon $\Delta$ of which a part of the boundary of $\Delta^{(1)}$ is as in Figure \ref{bijectionAB} (the $(i,j)$-coordinates are translated a bit). 
\begin{figure}[h!]
  \centering
    \includegraphics[height=3.5cm]{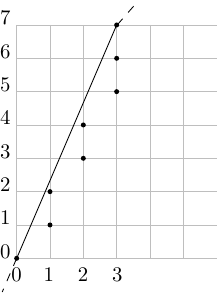}
		\caption{part of $\partial\Delta^{(1)}$}
		\label{bijectionAB}
\end{figure}

\noindent For this horizontal slice of the polygon, 
$$A^{(-)}=\{(0,2),(0,7),(2,4),(2,7),(4,6),(4,7)\} \quad \text{ and } \quad B^{(-)}=\{1,2,3,4,5,6\}.$$ 
The map $\alpha^{(-)}$ is defined as follows:
$$(0,2)\mapsto 1\ ,\ (0,7)\mapsto 2\ ,\ (2,4)\mapsto 3\ ,\ (2,7)\mapsto 4\ ,\ (4,6)\mapsto 5\ ,\ (4,7)\mapsto 6.$$
\end{example}

One can show that the $D_{j_1,j_2}$'s in Theorem \ref{thm_sfBn_toric} can be used to resolve $\mathcal{O}_{T'}$ as an $\mathcal{O}_{\mathbb{P}(\mathcal{E})}$-module, following Schreyer \cite{Schreyer}. For this one needs that the fibers of $\pi|_{T'}:T'\to \mathbb{P}^1$ have constant Betti numbers and that the corresponding resolutions are pure, but this can be verified. So it is justified to call the $B_{j_1,j_2}$'s the \emph{first scrollar Betti numbers of the toric surface} $\Tor(\Delta^{(1)})$, even though we will not push this discussion further.

\subsection{First scrollar Betti numbers of non-degenerate curves relative to the toric surface} \label{subsec2.4}

We will use the same set-up and assumptions as in the beginning of Section \ref{subsec2.3}. The assumption (i) implies that $\Delta^{(2)}\neq \emptyset$. 

\begin{definition} \label{def_P2}                                                                   
We say that $\Delta$ satisfies condition $\mathcal{P}_2(v)$ if for each lattice point $(i,j)$ of $\Delta$ and each horizontal line $L$, there exist two (not necessarily distinct) horizontal lines $M_1,M_2$, such that for all $w \in L \cap \Delta^{(2)} \cap \mathbb{Z}^2$
\[ \exists u_{i,j} \in M_1 \cap \Delta^{(1)} \cap \mathbb{Z}^2, \ v_{i,j} \in M_2 \cap \Delta^{(1)} \cap \mathbb{Z}^2 \ : \ (i,j) - w = (u_{i,j} - w) + (v_{i,j} - w ).\]
\end{definition}

\begin{remark} \label{rmk_equivdef} Write
$$i^{(--)}(j)=\min\{i\in\mathbb{Z}\,|\,(i,j)\in \Delta^{(2)}\} \text{ and } i^{(++)}(j)=\max\{i\in\mathbb{Z}\,|\,(i,j)\in \Delta^{(2)}\}$$
for all $j\in\{2,\ldots,\gamma-2\}$.
An equivalent definition is as follows: $\Delta$ satisfies condition $\mathcal{P}_2(v)$ if and only if for all $(i,j)\in\Delta$ and for all $j'\in\{2,\ldots,\gamma-2\}$, there exist $j_1,j_2\in\{1,\ldots,\gamma-1\}$ such that $j_1+j_2=j+j'$ and 
\begin{equation} \label{eq_P2}
i+[i^{(--)}(j'),i^{(++)}(j')]\subset [i^{(-)}(j_1),i^{(+)}(j_1)]+[i^{(-)}(j_2),i^{(+)}(j_2)].
\end{equation}

This condition is obviously satisfied for $(i,j) \in \Delta^{(1)}$ (take $j_1=j$ and $j_2=j'$). Moreover, the condition also holds if $(i,j)$ lies on the interior of a horizontal edge (i.e.\ the top or bottom edge) of $\Delta$. Indeed, assume for instance that $(i,j)$ lies in the interior of the top edge 
$[(i^-,j),(i^+,j)]$ of $\Delta$. We have that $$i^{(-)}(j'+1)+i^{(-)}(j-1)\leq i^{(--)}(j')+i^-+1\leq i^{(--)}(j')+i.$$ Hereby, the first inequality follows by replacing $L$ in the proof of Lemma \ref{lem_interior} by the half-closed line segment $[(i^{(--)}(j'),j'),(i^-,j)[$. Analogously, we get that $i^{(+)}(j'+1)+i^{(+)}(j-1)\geq i^{(++)}(j')+i$, so \eqref{eq_P2} follows for $j_1=j'+1$ and $j_2=j-1$.
\end{remark}

Although at first sight the condition $\mathcal{P}_2(v)$ might seem strong, it is not so easy to cook up counterexamples. However, there do exist instances of lattice polygons $\Delta$ for which the condition is not satisfied. The smallest example we have found is a polygon with $46$ interior lattice points and lattice width $10$.   

\begin{example} \label{exCEaligned}
Let $\Delta$ be as in Figure \ref{FigCEaligned} (the dashed line indicates $\Delta^{(1)}$).
\begin{figure}[h]
  \centering
    \includegraphics[height=4.5cm]{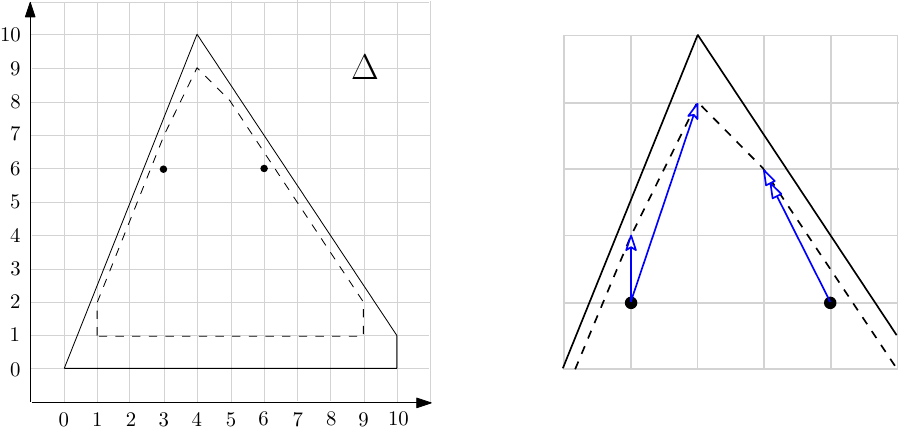} 
		\caption{A lattice polygon $\Delta$ that does not satisfy condition $\mathcal{P}_2(v)$}
		\label{FigCEaligned}
\end{figure}

\noindent We claim that $\Delta$ does not satisfy condition $\mathcal{P}_2(v)$. Indeed, take the top vertex $(i,j)=(4,10)$ of $\Delta$ and the horizontal line $L$ at height $6$. For the point $w\in L \cap \Delta^{(2)} \cap \mathbb{Z}^2$, consider the bold-marked lattice points $(3,6)$ and $(6,6)$ on $L$. In both cases, there is a unique decomposition of $(i,j) - w$:
\[ (1,4) = (0,1) + (1,3) \qquad \text{resp.} \qquad (-2,4) = (-1,2) + (-1,2).\]
So one sees that it is impossible to take the $u_{i,j}$'s and/or the $v_{i,j}$'s on the same line, which proves the claim.
\end{example}

\begin{theorem} \label{thm_tfBn}
If $\Delta$ satisfies condition $\mathcal{P}_2(v)$, there exist $\gamma-3$ effective divisors $D_\ell$ on $\mathbb{P}(\mathcal{E})$ (with $2\leq \ell\leq \gamma-2$) such that
\begin{itemize}
\item $C'$ is the (scheme theoretic) intersection of $T'$ and the divisors $D_\ell$,
\item $D_\ell\sim 2H-B_\ell R$ for all $\ell$, where $$B_\ell=i^{(++)}(\ell)-i^{(--)}(\ell)=-1+\sharp\{(i,j)\in \Delta^{(2)}\cap\mathbb{Z}^2\,|\,j=\ell\},$$
so $$\sum_{2\leq \ell\leq \gamma-2}\,B_\ell=\sharp(\Delta^{(2)}\cap\mathbb{Z}^2)-(\gamma-3).$$
\end{itemize}
\end{theorem}

\begin{proof}
The formula for the sum $\sum_\ell\,B_\ell$ is easily verified, so we focus on the other assertions.
Write
\[ f = \sum_{(i,j) \in \Delta \cap \mathbb{Z}^2} c_{i,j} x^iy^j \in k[x^{\pm 1},y^{\pm 1}] \]
and let $2\leq \ell\leq \gamma-2$ be an integer.

Take $\lambda\geq 0$ so that $S'=S(e_1+\lambda,\ldots,e_{\gamma-1}+\lambda)$ is smooth (and isomorphic to $\mathbb{P}(\mathcal{E})$) and define $\Delta'=\Delta[\lambda]$. We are going to use the inclusions 
$$C'\subset T'\subset S'\subset \mathbb{P}^{g-1+\lambda(\gamma-1)}$$ from \eqref{incl2}, where $T'=\Tor(\Delta'^{(1)})$. Write $(X_{i,j})_{(i,j)\in \Delta'^{(1)}\cap\mathbb{Z}^2}$ for the projective coordinates on $\mathbb{P}^{g-1+\lambda(\gamma-1)}$. 

Let $\ell\in\{2,\ldots,\gamma-2\}$ and denote the lattice points of $\Delta[2\lambda]^{(2)}$ on height $\ell$ by 
$w_0,\ldots,w_{B_\ell+2\lambda}$. If $w\in \{w_0,\ldots,w_{B_\ell+2\lambda}\}$ and $(i,j)\in\Delta$, then we claim that we can find  
$u_{i,j}, v_{i,j}\in \Delta'^{(1)}$ such that $(i,j)-w=(u_{i,j}-w)+(v_{i,j}-w)$, in such a way that their second coordinates are independent from $w$. 
Indeed, since $\Delta$ satisfies condition $\mathcal{P}_2(v)$, there exist $j_1,j_2\in\{1,\ldots,\gamma-1\}$ such that $j_1+j_2=j+\ell$ and $$i+[i^{(--)}(\ell),i^{(++)}(\ell)]\subset [i^{(-)}(j_1),i^{(+)}(j_1)]+[i^{(-)}(j_2),i^{(+)}(j_2)],$$ hence $$i+[i^{(--)}(\ell),i^{(++)}(\ell)+2\lambda]\subset [i^{(-)}(j_1),i^{(+)}(j_1)+\lambda]+[i^{(-)}(j_2),i^{(+)}(j_2)+\lambda].$$ This implies that we can take $u_{i,j}$ and $v_{i,j}$ with second coordinates $j_1$ and $j_2$. 
Define $$Q_w=\sum_{(i,j) \in \Delta \cap \mathbb{Z}^2} c_{i,j} X_{u_{i,j}} X_{v_{i,j}} \in k[X_{i,j}]_{(i,j)\in \Delta'^{(1)}\cap\mathbb{Z}^2}.$$ 
A consequence of the choice of $u_{i,j}, v_{i,j}$ is that $$X_{w_{s}}Q_{w_{r+1}}-X_{w_{s+1}}Q_{w_r} \in \mathcal{I}(S')$$ (rather than just $\mathcal{I}(T')$) for all 
$r\in \{0,\ldots,B_\ell+2\lambda-1\}$ and $s\in \{0,\ldots,B_\ell+\lambda-1\}$. Since 
$$\frac{X_{w_1}}{X_{w_0}} = \frac{X_{w_2}}{X_{w_1}} = \hdots = \frac{X_{w_{B_\ell+\lambda}}}{X_{w_{B_\ell+\lambda - 1}}}$$ is a local parameter for the $(\gamma-2)$-plane $R_{(0:1)} = \pi^{-1}(0:1) \subset S'$, it follows that the $R_{(0:1)}$-orders of
\begin{equation} \label{ZQws}
\mathcal{Z}(Q_{w_1}), \mathcal{Z}(Q_{w_2}), \dots, \mathcal{Z}(Q_{w_{B_\ell+2\lambda}})
\end{equation}
increase by $1$ at each step. For a similar reason, with $R_{(1:0)} = \pi^{-1}(1:0)\subset S'$, the
$R_{(1:0)}$-orders of \eqref{ZQws} decrease by $1$ at each step. We conclude that there exists an effective divisor $D_\ell$ such that for all 
$i\in\{0,\dots,B_\ell+2\lambda\}$ we have
\begin{equation} \label{ZQwdecomp}
\mathcal{Z}(Q_{w_i}) \, = \, i \cdot R_{(0:1)} \, + \, (B_\ell + 2\lambda - i) \cdot R_{(1:0)} \, + \, D_{\ell}
\end{equation}
on $S'$. The divisor $D_\ell$ is in fact the divisor of $S'$ cut out by the quadrics in \eqref{ZQws}. Using \eqref{ZQwdecomp} and Remark \ref{rmk_ezero}, we get that 
$$D_\ell \sim 2H' - (B_\ell+2\lambda) R=2H-B_\ell R,$$ 
so it is sufficient to show that the quadrics $Q_w$ (where $w$ ranges over 
$\Delta[2\lambda]^{(2)}\cap\mathbb{Z}^2$) cut out $C'$ from $T'$. If $\lambda=0$, this follows from \cite[Theorem 3.3]{CaCo2}.

Before we prove this, we need to introduce one more notion: for each lattice polygon $\Gamma$ with $\Gamma^{(1)}$ two-dimensional, write $\Gamma^{max}$ to denote the largest lattice polygon with interior lattice polygon equal to $\Gamma^{(1)}$, so $\Gamma^{max}\supset \Gamma$. If $\Gamma^{(1)}$ is the intersection of the half-planes $$H_t=\{P\in\mathbb{R}^2\,|\,\langle P,v_t\rangle\geq a_t\}$$ (where $\langle\cdot,\cdot\rangle$ denotes the standard inner product, $v_t$ is a primitive inward pointing normal vector of an edge of $\Gamma^{(1)}$ and $a_t\in\mathbb{Z}$), then $$\Gamma^{max}=\cap_t\,H_t^{(-1)} \quad \text{with} \quad H_t^{(-1)}=\{P\in\mathbb{R}^2\,|\,\langle P,v_t\rangle\geq a_t-1\}.$$ We will use the following two properties (see \cite[Section 2]{CaCo1} for other properties of $\Gamma^{max}$): 
\begin{itemize}
\item If $\Gamma^{(2)}\neq \emptyset$, then $2\Gamma^{(1)}=\Gamma^{(2)}+\Gamma^{max}$, since both lattice polygons are defined by the half-planes $2H_t=\{P\in\mathbb{R}^2\,|\,\langle P,v_t\rangle\geq -2a_t\}$. 
\item If $\Phi_1,\Phi_2$ are lattice polygons such that $\Phi_1 + \Gamma \subset \Phi_2 + \Gamma^{max}$, then $\Phi_1\subset \Phi_2$ if $\Phi_2$ satisfies the following condition: it is the intersection of half-planes $H'_t$ with $H'_t$ of the form $$\{P\in\mathbb{R}^2\,|\,\langle P,v_t\rangle\geq b_t\}$$ for some $b_t\in\mathbb{Z}$ (hence parallel to $H_t$).  

Indeed, if $\Phi_1\not\subset \Phi_2$, take a lattice point $P$ in $\Phi_1\setminus\Phi_2$. Then $\langle P,v_t\rangle<b_t$ for some value of $t$. 
Take $Q\in\Gamma$ with $\langle Q,v_t\rangle=a_t-1$ (this is always possible). We have that $P+Q\in \Phi_1+\Gamma$ and $\langle P+Q,v_t\rangle<a_t+b_t-1$, but $\Phi_2 + \Gamma^{max}$ is the intersection of the half-planes $H''_t=\{P\in\mathbb{R}^2\,|\,\langle P,v_t\rangle\geq a_t+b_t-1\}$, so $P+Q\not\in \Phi_2 + \Gamma^{max}$, a contradiction.
\end{itemize}

Now take $F\in \mathcal{I}(C')$ homogeneous of degree $d$. If $$\xi : k[X_{i,j}]_{(i,j)\in \Delta'^{(1)}\cap\mathbb{Z}^2} \to k[x^{\pm 1},y^{\pm 1}]$$ is the ring morphism that maps $X_{i,j}$ to $x^i y^j$, then $\xi(F)$ has to be of the form $cf$ for some $c\in k[x^{\pm 1},y^{\pm 1}]$. The Newton polygon of $cf$ is equal to $\Delta(c)+\Delta$, while the Newton polygon $\Delta(\xi(F))$ is contained in 
$$d\Delta'^{(1)}=(d-2)\Delta'^{(1)}+\Delta'^{(2)}+\Delta'^{max}=(d-2)\Delta'^{(1)}+\Delta[2\lambda]^{(2)}+\Delta^{max}$$ 
(here, we use the first property of maximal polygons with $\Gamma=\Delta'$). So we obtain that 
$$\Delta(c)+\Delta\subset (d-2)\Delta'^{(1)}+\Delta[2\lambda]^{(2)}+\Delta^{max}.$$ Now we can use the second property of maximal polygons with $\Phi_1=\Delta(c)$, $\Phi_2=(d-2)\Delta'^{(1)}+\Delta[2\lambda]^{(2)}$ and $\Gamma=\Delta$. Note that $\Phi_2$ might have a horizontal (top or bottom) edge while $\Delta^{(1)}$ has not, but this is not an issue (since $\Delta^{(1)}\not\cong (d-3)\Sigma$). It follows that 
$$\Delta(c)\subset (d-2)\Delta'^{(1)}+\Delta[2\lambda]^{(2)}.$$
So we can write $$c=\sum_{w=(i,j)\in \Delta[2\lambda]^{(2)}\cap\mathbb{Z}^2}\, g_{i,j}x^i y^j$$ for polynomials $g_{i,j}\in k[x,y]$ with 
$\Delta(g_{i,j})\subset (d-2)\Delta'^{(1)}$. For all lattice points $w=(i,j)\in \Delta[2\lambda]^{(2)}\cap\mathbb{Z}^2$, there is a homogeneous polynomial $G_{i,j}\in k[X_{i,j}]_{(i,j)\in \Delta'^{(1)}\cap\mathbb{Z}^2}$ such that $\xi(G_{i,j})=g_{i,j}$. On the other hand, $\xi(Q_w)=x^i y^j f$, hence $$\xi(F)=cf=\sum_{w=(i,j)\in\Delta[2\lambda]^{(2)}\cap\mathbb{Z}^2}\, \xi(G_{i,j})\xi(Q_w).$$ So $F-\sum_{w=(i,j)}\, G_{i,j}Q_w$ belongs to the kernel of the map $\xi$, which implies that it is contained in $\mathcal{I}_d(T')$, which is what we wanted to prove. 
\end{proof}

\subsection{First scrollar Betti numbers for non-degenerate curves}

We are ready to prove the main result of this section, by combining the results from Sections \ref{subsec2.3} and \ref{subsec2.4}. 

\begin{theorem} \label{thm_main}
Let $\Delta$ be a lattice polygon with $\lw(\Delta)\geq 4$ such that $\Delta^{(1)}\not\cong \Upsilon$ and $\Delta^{(1)}\not\cong (d-3)\Sigma$ for any integer $d\geq 3$. Moreover, assume that $\Delta$ satisfies the conditions $\mathcal{P}_1(v)$ and $\mathcal{P}_2(v)$, where $v$ is a lattice width direction. Let $C$ be a $\Delta$-non-degenerate curve and let $g_\gamma^1$ be the combinatorial gonality pencil on $C$ corresponding to $v$ (with $\gamma=\lw(\Delta)$). 
Then the first scrollar Betti numbers of $C$ with respect to $g_\gamma^1$ are given by 
$$\{B_\ell\}_{\ell\in\{2,\ldots,\gamma-2\}}\cup {\{B_{j_1,j_2}\}}_{\substack{j_1,j_2\in\{1,\ldots,\gamma-1\} \\ j_2-j_1\geq 2}}.$$
\end{theorem}

\begin{proof}
We use the notations and set-up from \eqref{incl2}. Theorem \ref{thm_sfBn_toric} and Theorem \ref{thm_tfBn} imply that there exist divisors $D_\ell\sim 2H-B_\ell R$ on $\mathbb{P}(\mathcal{E})$, with $\ell\in\{2,\ldots,\gamma-2\}$, and divisors $D_{j_1,j_2}\sim 2H-B_{j_1,j_2}R$ on $\mathbb{P}(\mathcal{E})$, with $j_1,j_2\in\{1,\ldots,\gamma-1\}$ and $j_2-j_1\geq 2$, such that $C'$ is the scheme-theoretic intersection of these divisors. Since moreover 
$$\sum_{\ell\in\{2,\ldots,\gamma-2\}}\,B_\ell\ + \sum_{\substack{j_1,j_2\in\{1,\ldots,\gamma-1\} \\ j_2-j_1\geq 2}}\,B_{j_1,j_2}
=(\gamma-3)g-(\gamma^2-2\gamma-3),$$
we can use Theorem \ref{thm_fsBn} to conclude the proof. 
\end{proof}

We believe that the above theorem is of independent interest. For instance it is not well-understood which sets of (first) scrollar Betti numbers are possible for canonical curves of a given genus and gonality, and our result can be used to prove certain existence results. It has been conjectured that ``most'' curves have so-called balanced (first) scrollar Betti numbers, meaning that $\max |b_i-b_j|\leq 1$, see \cite{BoppHoff} and the references therein. Non-degenerate curves are typically highly non-balanced, since one expects the $B_{j_1,j_2}$'s to be about twice the $B_\ell$'s. 

\begin{example}
Consider the following lattice polygons $\Delta_1$ and $\Delta_2$ of lattice width $7$ (and lattice width direction $v=(1,0)$), which only differ from each other at the right hand side. 
\begin{center}
\includegraphics[height=4.5cm]{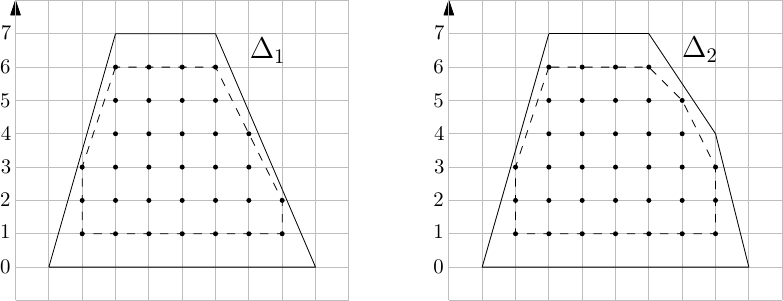}
\end{center}

The polygon $\Delta_1$ does not satisfy all the combinatorial constraints of Theorem \ref{thm_main}, since condition $\mathcal{P}_2(v)$ does not hold: $\mathcal{S}_{2,6}^{(-)}=\{1\}$ and $\mathcal{S}_{2,6}^{(+)}=\{2\}$. Although we have not pursued this, we believe that the conditions $\mathcal{P}_1(v)$ and $\mathcal{P}_2(v)$ are always fulfilled if $\gamma<7$. 

On the other hand, the polygon $\Delta_2$ meets all the conditions from the statement, and so we can apply the theorem. The first scrollar Betti numbers of a $\Delta_2$-degenerate curve are as follows: 
$$\begin{array}{lllll}
B_{1,3}=12 & B_{1,4}=10 & B_{1,5}=10 & B_{1,6}=9 & B_{2,4}=10 \\ 
B_{2,5}=10 & B_{2,6}=9 & B_{3,5}=8 & B_{3,6}=8 & B_{4,6}=7 \\
B_2=4 & B_3=4 & B_4=4 & B_5=3 & 
\end{array}$$
The sum of these numbers is $108$, which agrees with $(\gamma-3)g-(\gamma^2-2\gamma-3)$ for $g=35$ and $\gamma=7$.
\end{example}

\section{Intrinsicness on \texorpdfstring{$\mathbb{P}^1\times \mathbb{P}^1$}{P1*P1}} \label{sec_intrinsicness}

\begin{theorem} \label{P1xP1theorem}
   Let $f \in k[x^{\pm 1},y^{\pm 1}]$ be non-degenerate with respect to its (two-dimensional) Newton polygon
 $\Delta = \Delta(f)$, and assume that $\Delta \not \cong 2\Upsilon$. Then $U(f)$ is birationally equivalent to a smooth projective genus $g$ curve in $\mathbb{P}^1 \times \mathbb{P}^1$
 if and only if $\Delta^{(1)} = \emptyset$ or $\Delta^{(1)} \cong \square_{a,b}$ for some integers $a\geq b \geq 0$, necessarily satisfying $g=(a+1)(b+1)$.
\end{theorem}

\begin{proof}
We may assume that $U(f)$ is neither rational, nor elliptic or hyperelliptic
(and hence that $\Delta^{(1)}$ is two-dimensional) because such curves admit smooth complete models in $\mathbb{P}^1 \times \mathbb{P}^1$.
So for the `if' part we can assume that $b \geq 1$. But then $\text{Tor}(\Delta^{(1)}) \cong \mathbb{P}^1 \times \mathbb{P}^1$, and
the statement follows using the canonical embedding \eqref{canemb}.

The real deal is the `only if' part.
At least, if a curve $C/k$ is birationally equivalent to a (non-rational, non-elliptic, non-hyperelliptic)
smooth projective curve in $\mathbb{P}^1 \times \mathbb{P}^1$, then
it is $\Delta'$-non-degenerate with $\Delta' = [-1, a+1] \times [-1, b+1]$ for $a \geq b \geq 1$: this
follows from the material in \cite[Section 4]{CaCo1} (one can use an automorphism of $\mathbb{P}^1 \times
\mathbb{P}^1$ to ensure appropriate behavior with respect to the toric boundary). Note that
$\Delta'^{(1)} = \square_{a,b}$. The geometric
genus of $C$ equals $g=(a+1)(b+1)$ by \cite{Khov} and its gonality equals $\gamma = b + 2$ by \cite[Cor.\ 6.2]{CaCo1}.
We observe that
\begin{itemize}
\item $g$ is composite,
\item $C$ has Clifford dimension equal to $1$ by \cite[Theorem 8.1]{CaCo1},
\item the scrollar invariants of $C$ (with respect to any gonality pencil) are all equal to $g/(\gamma - 1) - 1$; indeed,
by \cite[Theorem 6.1]{CaCo1} every gonality pencil is computed
by projecting along some lattice width direction $v$; if $a > b$ then the only pair of lattice width directions
is $\pm (1, 0)$; from \cite[Theorem 9.1]{CaCo1} we find
that the corresponding scrollar invariants
are $a, a, \dots, a$; if $a = b$ we also have the pair $\pm (0,1)$, giving rise to the same scrollar invariants,
\item if $\gamma \geq 4$ then the first scrollar Betti numbers (with respect to any gonality pencil) take
exactly two
distinct values: $2g/(\gamma - 1) - 2$ and $g/(\gamma - 1) - 3$. Indeed, $\Delta'$ satisfies condition $\mathcal{P}_1(v)$ (see Example \ref{easyexample}), but also condition $\mathcal{P}_2(v)$: 
take $(j_1,j_2)=(j,\ell)$ if $j\in\{-1,\ldots,b+1\}$, $(j_1,j_2)=(j+1,\ell-1)$ if $j=-1$ and $(j_1,j_2)=(j-1,\ell+1)$ if $j=b+1$. By Theorem \ref{thm_main}
we find that these numbers are $2a, 2a, \dots, 2a, a-2,a-2, \dots, a-2$.
\end{itemize}
A first consequence is that $U(f)$ admits a combinatorial gonality pencil. Indeed, $\Delta$ cannot
be of the form $2\Upsilon$ (excluded in the statement of the theorem), nor of the form $d\Sigma$ for some $d \geq 2$:
the cases $d =2$ and $d=3$
correspond to rational and elliptic curves (excluded at the beginning of this proof), the case $d=4$ corresponds to curves of genus $3$ (not composite),
and the cases where $d \geq 5$ correspond to curves of Clifford dimension $2$.

Without loss
of generality we may then assume that $v=(1,0)$ and $\Delta \subset \{ \, (i,j)\in \mathbb{R}^2\,|\, 0\leq j\leq \gamma \, \}$, so
that our gonality pencil corresponds to vertical projection.
By \cite[Theorem 9.1]{CaCo1}, the numbers
$E_\ell = -1+\sharp\{(i,j)\in \Delta^{(1)}\cap\mathbb{Z}^2\,|\,j=\ell\}$
(for $\ell = 1, \dots, \gamma - 1$) are the corresponding scrollar invariants. Hence the $E_\ell$'s must all be equal
to $E:=g/(\gamma - 1) - 1\geq 1$.

This already puts severe restrictions on the possible shapes of $\Delta^{(1)}$.
By horizontally shifting and skewing
we may assume that the lattice points at height $j=1$ are $(0,1), \dots, (E, 1)$ and that
the lattice points at height $j=2$ are $(0,2), \dots, (E,2)$. If $\gamma = 3$, it follows
that $\Delta^{(1)}$ has the desired rectangular shape, so we may suppose that $\gamma \geq 4$.
Then by horizontally flipping if needed, we can assume that the lattice points at height $j=3$
are $(i,3), \dots, (E + i, 3)$ for some $i \geq 0$. Now $i \geq 2$ is impossible, for this
would introduce a new lattice point at height $j=2$; thus $i = 0$ or $i=1$.
Continuing this type of reasoning, we obtain that the lattice points of $\Delta^{(1)}$
can be seen as a pile of $n$ blocks of respectively $m_1, \dots, m_n$ sheets, where each block is shifted to the right
over a distance $1$ when compared to its predecessor.

\begin{figure}[h]
  \centering
    \includegraphics[height=5cm]{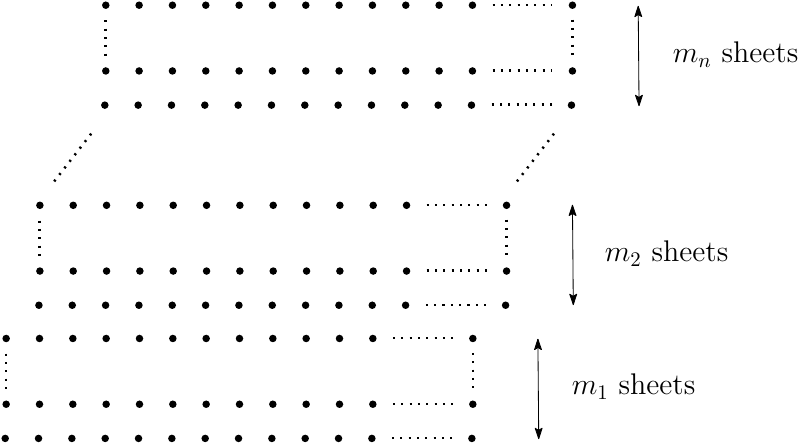} 
		\caption{lattice points of $\Delta^{(1)}$ in sheets}
		\label{sheets}
\end{figure}

We need to show that $n = 1$, because then $\Delta^{(1)}$ has the desired rectangular shape. We will first prove that $\Delta$ statisfies condition $\mathcal{P}_1(v)$. Since we have that $i^{(+)}(j)-i^{(-)}(j)=E$ for each value of $j$, the inequality $$\epsilon^{(-)}_{j_1,j_2,r}+\epsilon^{(+)}_{j_1,j_2,r}\leq 1$$ holds (so never $\epsilon^{(-)}_{j_1,j_2,r}=\epsilon^{(+)}_{j_1,j_2,r}=1$) for all $j_1,j_2\in\{2,\ldots,\gamma-2\}$ with $j_2-j_1\geq 2$ and $r\in \left\{1,\ldots,\left\lfloor \frac{j_2-j_1}{2}\right\rfloor\right\}$. This implies that $$\mathcal{S}^{(-)}_{j_1,j_2}\cup \mathcal{S}^{(+)}_{j_1,j_2}=\left\{1,\ldots,\left\lfloor \frac{j_2-j_1}{2}\right\rfloor\right\}.$$ Now assume that $\mathcal{S}^{(-)}_{j_1,j_2}$ and $\mathcal{S}^{(+)}_{j_1,j_2}$ are non-empty and disjoint. In this case, we can take $r,s\in \left\{1,\ldots,\left\lfloor \frac{j_2-j_1}{2}\right\rfloor\right\}$ such that $$\epsilon^{(-)}_{j_1,j_2,r}=\epsilon^{(+)}_{j_1,j_2,s}=0 \quad \text{ and } \quad \epsilon^{(+)}_{j_1,j_2,r}=\epsilon^{(-)}_{j_1,j_2,s}=1.$$ If $r<s$, we get that $$i^{(-)}(j_1 + s) + i^{(-)}(j_2 - s) > i^{(-)}(j_1) + i^{(-)}(j_2)$$ and $$i^{(+)}(j_1) + i^{(+)}(j_2) > i^{(+)}(j_1 + r) + i^{(+)}(j_2 - r).$$
Subtracting $E$ from both sides of the latter equation yields
$$i^{(-)}(j_1) + i^{(-)}(j_2) > i^{(-)}(j_1 + r) + i^{(-)}(j_2 - r),$$
so $$i^{(-)}((j_1+r)+(s-r))+i^{(-)}((j_2-r)-(s-r))\geq i^{(-)}(j_1+r)+i^{(-)}(j_2-r)+2,$$ which is in contradiction with Lemma \ref{lem_interior}. A similar contradiction can be obtained if $r>s$. 

Now let's prove that $\Delta$ also satisfies property $\mathcal{P}_2(v)$, where we assume that $n\geq 2$. By Remark \ref{rmk_equivdef} and a symmetry consideration (rotation over $180^\circ$), it suffices to check the condition for lattice points $(i,j)$ that lie on the left side of the boundary of $\Delta$ (and even of $\Delta^{max}$). Take $\ell\in\{2,\ldots,\gamma-2\}$, 
$w=(i^{(--)}(\ell),\ell)$ and 
$u_{i,j}=(i_1,j_1),v_{i,j}=(i_2,j_2)\in\Delta^{(1)}$ such that $(i,j)-w=(u_{i,j}-w)+(v_{i,j}-w)$, hence $j_1+j_2=j+\ell$. It is sufficient to prove that 
\begin{equation} \label{eqP2}
i+i^{(++)}(\ell)\leq i^{(+)}(j_1)+i^{(+)}(j_2). 
\end{equation}                          

First assume that $|j-\ell|>|j_2-j_1|$. If 
$j\in\{1,\ldots,\gamma-1\}$, then 
\begin{eqnarray*}
(i+E)+i^{(++)}(\ell) &\leq& (i^{(-)}(j)+E)+i^{(+)}(\ell) \\ &=& i^{(+)}(j)+i^{(+)}(\ell) \\
&\leq& i^{(+)}(j_1)+i^{(+)}(j_2)+1,
\end{eqnarray*}
where we use Lemma \ref{lem_interior} for the last inequality. Since $E\geq 1$, the desired inequality \eqref{eqP2} follows. We still need to check \eqref{eqP2} for points $(i,j)\in\partial\Delta$ with $j=0$ and $j=\gamma$, in particular $i=-1$ resp.\ $i=n-1$ because we can assume that $(i,j)$ lies on the left side of the boundary of $\Delta^{max}$. 
\begin{itemize}
\item If $(i,j)=(-1,0)$, the line segment $L$ between $(i+E,j)=(E-1,0)$ and $w'=(i^{(++)}(\ell),\ell)\in \Delta^{(2)}$ intersects the horizontal lines on heights $j_1$ and $j_2$ in points that belong to $\Delta^{(1)}$. Using a similar argument as in the proof of Lemma 
\ref{lem_interior}, we obtain that $(i+E)+i^{(++)}(\ell)\leq i^{(+)}(j_1)+i^{(+)}(j_2)+1$, which gives us \eqref{eqP2} using $E\geq 1$.       
\item Analogously, we can handle the case $(i,j)=(n-1,\gamma)$: the line segment $L$ between $(i+E,j)=(n+E-1,\gamma)$ and $w'$ will intersect the horizontal lines on heights $j_1$ and $j_2$ in points that are contained in $\Delta^{(1)}$ and \eqref{eqP2} follows.    
\end{itemize} 

If $|j-\ell|=|j_2-j_1|$, we may assume that $j_1=j\in\{1,\ldots,\gamma-1\}$ and $j_2=\ell\in\{2,\ldots,\gamma-2\}$. But then 
the inequalities $i\leq i^{(-)}(j_1)\leq i^{(+)}(j_1)$ and $i^{(++)}(\ell)\leq i^{(+)}(j_2)$ yield \eqref{eqP2}.

We still have to consider the case where $|j-\ell|<|j_2-j_1|$, which implies that $j\in\{1,\ldots,\gamma-1\}$. By Lemma \ref{lem_interior}, we have that $$i^{(-)}(j)+i^{(-)}(\ell)\leq i^{(-)}(j_1)+i^{(-)}(j_2)+1,$$ hence
\begin{eqnarray*}
i+i^{(++)}(\ell) &\leq& i^{(-)}(j)+i^{(+)}(\ell) \\ &=& i^{(-)}(j)+i^{(-)}(\ell)+E \\ 
&\leq& i^{(-)}(j_1)+i^{(-)}(j_2)+E+1 \\ &\leq& i^{(+)}(j_1)+i^{(+)}(j_2).
\end{eqnarray*}

Since both the conditions $\mathcal{P}_1(v)$ and $\mathcal{P}_2(v)$ hold for $\Delta$, we can apply Theorem \ref{thm_main}. If $n \geq 2$, then there
is at least one first scrollar Betti number having value $E - 1 = g/(\gamma - 1) - 2$, for instance 
$B_2$. This is distinct from both $2g/(\gamma - 1) - 2$ and $g/(\gamma - 1) - 3$: contradiction. Therefore $n=1$, i.e.\ $\Delta^{(1)}$ has
the requested rectangular shape.
\end{proof}

\section{Open questions}

Here are two interesting open questions related to this paper:
\begin{enumerate}
\item In Section \ref{sec2}, we gave a combinatorial interpretation for the first scrollar Betti numbers of $\Delta$-non-degenerate curves $C$ in terms of the combinatorics of $\Delta$, in case $\Delta$ satisfies the condition $\mathcal{P}_1(v)$ (see Definition \ref{def_P1}) and $\mathcal{P}_2(v)$ (see Definition \ref{def_P2}). Can this be generalized to all polygons $\Delta$? 
There seems to be no geometric reason why this wouldn't be the case, but we did not succeed to get rid of the condition. 
\item In Theorem \ref{P1xP1theorem} of Section \ref{sec_intrinsicness}, we showed that non-degenerate curves on $\mathbb{P}^1\times \mathbb{P}^1$ have an intrinsic Newton polygon (at least, if $g\neq 4$). Can this be generalized to $\Delta$-non-degenerate curves on Hirzebruch surfaces $\mathcal{H}_n$? In this case, we expect $\Delta^{(1)} = \emptyset$ or $\Delta^{(1)} \cong \conv\{(0,0),(a+nb,0),(a,b),(0,b)\}$ for some integers $a,b,n \geq 0$ .  
\end{enumerate}

\vspace{5mm}
\noindent \textsc{Vakgroep Wiskunde, Universiteit Gent}\\
\noindent \textsc{Krijgslaan 281, 9000 Gent, Belgium}\\
\vspace{-0.4cm}

\noindent \textsc{Departement Elektrotechniek, KU Leuven and iMinds}\\
\noindent \textsc{Kasteelpark Arenberg 10/2452, 3001 Leuven, Belgium}\\
\vspace{-0.4cm}

\noindent \emph{E-mail address:} \verb"wouter.castryck@gmail.com"\\

\noindent \textsc{Departement Wiskunde, KU Leuven}\\
\noindent \textsc{Celestijnenlaan 200B, 3001 Leuven, Belgium}\\
\vspace{-0.4cm}

\noindent \emph{E-mail address:} \verb"filip.cools@wis.kuleuven.be"\\

\end{document}